\documentclass{elsart}

\newtheorem{theorem}{Theorem}

\usepackage{amssymb}

\begin{document}

\begin{frontmatter}



\title{LAD Regression and Nonparametric Methods for Detecting 
Outliers and Leverage Points}


\author{Susana Faria${}^1$ and Giuseppe Melfi${}^2$}

\address{${}^1$University of Minho, DMCT\\
P-4800-058 {Guimar$\tilde{a}$es}, Portugal\\
sfaria@mct.uminho.pt\\
\mbox{   }\\
${}^2$
Universit\'e de Neuch\^atel, Institut de Statistique\\
Espace de l'Europe 4, CH--2002 Neuch\^atel, Switzerland\\
Giuseppe.Melfi@unine.ch}

\begin{abstract}
The detection of influential observations for the
standard least squares regression model is a question that has
been extensively studied.  LAD regression
diagnostics offers alternative approaches whose main feature
is the robustness.
In this paper a new approach for nonparametric
detection of influencial observations in LAD regression
models is presented and compared with other classical methods of
diagnostics.
\end{abstract}

\begin{keyword}
Least Absolute Deviations Regression \sep Robustness \sep
Outliers \sep Leverage Points.
\end{keyword}
\end{frontmatter}


\section{Introduction}\label{sectionone}
The robustness of LAD to low-leverage outliers, and its susceptibility
to high-leverage outliers has been extensively studied in 
literature~\cite{dodge1,dodge2,dodge3}.
In this paper we propose  a method for nonparametic detection of such 
influencial observations
by the use of a technique derived from LAD regression. Robust methods 
based on the $L_1$-norm have been proposed for example in
\cite{hadi,rousseeuwleroy}. The approach presented here considers
suitable perturbations of a given data set and allows a detection os
high-leverage observations and outliers from a new viewpoint. These methods
answer to natural requirements for robustness, and give a new tool
for the analysis of data.

Let $S\subset\mathbb R^{p+1}$ be a finite discrete set of points. In 
statistics, such a set may represent observations in 
$p+1$ variables.
Denote the elements
of $S$ as $(x_{i1},\dots,x_{ip},y_i)$, where the last variable is 
explained from the preceding ones by a linear regression model:
$$y_i=\beta_0+\sum_{j=1}^p\beta_jx_{ij}+\varepsilon_i\qquad 
\mbox{ for }i=1,\ldots,n,$$
where $p$ is the number of explanatory variables, 
$\varepsilon_i$ are error terms, 
or deviations, and $n$ is the number of observations.
The LAD regression model is determined by minimizing the 
sum of the absolute deviations, i.e., the vector 
$(\beta_0,\beta_1,\dots,\beta_p)\in\mathbb R^{p+1}$ 
is determined by minimizing on $\beta_0,\beta_1,\dots,\beta_p$ the function
\begin{equation}\label{min}
F(\beta_0,\beta_1,\dots,\beta_p):=
\sum_{i=1}^n\left|y_i-\beta_0-\sum_{j=1}^p\beta_jx_{ij}\right|.
\end{equation}


When a linear LAD regression model is fitted, the hyperplane always 
passes through at least
$p+1$ points~\cite{AD}, 
although the solution may be non-unique. 

For our purposes, we assume
that for every dataset and for each subset of it we deal with,
the hyperplane which fits the 
linear LAD regression model is unique, as well as the observation with
maximal absolute deviation.
These assumptions are reasonable for datasets whose size is sufficiently large
and/or the data contain sufficient significant digits.
We suppose also 
that the dataset is such that every $p+2$ points are not in the same 
hyperplane. With these assumptions the linear LAD 
regression model is unique and it passes through exactly $p+1$ points. 
Furthermore, if $n>p+1$, there is always a point which does not 
belong to the regression hyperplane, so having a positive absolute
deviation.

Consider the $n$ datasets composed by all possible subsets of $S$ of 
size $n-1$. Under the above assumptions for each dataset we have a 
unique solution. For each case, we assign the score 1 to each point 
through which the fitted model passes
and 0 to the other points. We define the final score of each point as 
the sum of scores over all models fitting the $n$ datasets. 
This score is produced by the repeated use of the same points, 
each time considering a different subset  of the original data set, 
so in a certain sense by bootstrapping the linear LAD regression model.

The point $(x_{k1},x_{k2},\dots,x_{kp},y_k)$ will be also denoted by $k$
and  its score  will be denoted by $L(k)$.

Similarly, we may define another complementary score function, denoted by
$O(k)$, in the following way.
Consider again the $n$ datasets composed 
by all possible subsets of $S$ of size $n-1$.
For each subset we consider the LAD regression line and 
we give the score $1$ to the unique (according to the above assumptions)
point which maximizes the
absolute distance from the LAD regression line.

We define
the score $O(k)$ as the sum (over all $n$ possible subsets of $S$
of size $n-1$)
of scores arising from the LAD regression lines.  

In Section~\ref{remarks} we discuss some elementary properties 
of the LAD regression
model. These properties will justify our algorithms 
for the detection of outliers and leverage 
points presented in Section~\ref{algorithms}. 
In Section~\ref{examples}
we discuss some examples, and compare the 
results with those obtained using
other classical methods.

\section{Preliminary Considerations}\label{remarks}

Under the above assumptions, the sum 
of $L$ scores of all points is $n(p+1)$, and the sum of $O$ scores is $n$,
so, under the random variable viewpoint, $E(L(k))= p+1$
and $E(O(k))=1$. Now suppose that we have a set
of observations, all concentrated in a region and an isolated observation 
horizontally very far  from the others but such that the
line of the LAD regression model 
will pass through it (a typical leverage point). 
It is likely that the $L$ score of this observation will be quite high. 
So a large $L$ score is synonymous
of leverage point. On the other hand, suppose we have a dataset
in which all points 
are roughly in 
a hyperplane, and a further point far above them (an outlier). 
The LAD regression model will be very near to this hyperplane, 
and the score $L$ of the outlier 
will be probably zero, but the score $O$ will be probably $n-1$.

To justify these arguments we state the following theorems.

\begin{theorem}
Let $(x_{11},\dots,x_{1p},y_1),\dots,(x_{n1},\dots,x_{np},y_n)$ 
be $n$ points in $\mathbb R^{p+1}$.
Let $(x_{{n+1}1},\dots,x_{{n+1}p},y_{n+1})$ be an additional
point, such that $(x_{{n+1}1},\dots,x_{{n+1}p})$ belongs to the
interior of the convex hull determined by
the set $\{(x_{11},\dots,x_{1p})$,
$\dots$, $(x_{n1},\dots,x_{np})\}$.
If $x_{{n+1}1},\dots,x_{{n+1}p}$ are fixed and 
$|y_{n+1}|$ is sufficiently large, a hyperplane relative to a
linear LAD regression model does not 
passes
through  $(x_{{n+1}1},\dots,x_{{n+1}p},y_{n+1})$.
\end{theorem}

\noindent\textbf{Proof.~}
Suppose $p=1$ and that for $i=1,\dots n$, \ $c<y_i<d$. 
The convex hull hypothesis reduces to
$$a=\min_{i=1,\dots,n}\{x_i\}<x_{n+1}< \max_{i=1,\dots,n}\{x_i\}=b.$$
Let $y=\beta_0+\beta_1 x$ be the line of the linear LAD regression model.
Let $\ell_1$ be the horizontal line $y=(c+d)/2$. The sum of the absolute
deviations of $(x_{11},\dots,x_{1p},y_1),\dots,
(x_{{n+1}1},\dots,x_{{n+1}p},y_{n+1})$ does not exceed $(d-c)(n+1)+|y_{n+1}|$
$+|d|+|c|$.
So
$$F(\beta_0,\beta_1)<(d-c)(n+1)+|d|+|c|+|y_{n+1}|.$$
On the other hand, if a line $y=\beta_0^*+\beta_1^*x$ passes through 
$(x_{n+1},y_{n+1})$, if it is a linear LAD regression model, it will 
pass also through another point $(x_i,y_i)$ for a suitable $i$,
so for sufficiently large $|y_{n+1}|$,
$|\beta_1^*|>|y_{n+1}|/(b-a)$. Hence there exists $\alpha>1$
such that for sufficiently large $|y_{n+1}|$,
$$F( \beta_0^*,\beta_1^*)>\alpha|y_{n+1}|.$$
For sufficiently large $|y_{n+1}|$ we have
$$F(\beta_0^*,\beta_1^*)>\alpha|y_{n+1}|>(d-c)(n+1)+|d|+|c|+|y_{n+1}|>
F(\beta_0,\beta_1),$$
and therefore, the line relative to the the LAD regression model cannot pass
through $(x_{n+1},y_{n+1})$.


For $p>1$ the proof is similar.
$\Box$\medskip

\begin{theorem}
Let $(x_{11},\dots,x_{1p},y_1),\dots,(x_{n1},\dots,x_{np},y_n)$ 
be $n$ points in $\mathbb R^{p+1}$.
Let $(x_{{n+1}1},\dots,x_{{n+1}p},y_{n+1})$ be an additional
point.
If $y_{n+1}$ is fixed and $\sum_{i=1}^p|x_{n+1,i}|$ is sufficiently large,
a linear LAD regression 
model will pass through  $(x_{{n+1}1}$, $\dots$,
$x_{{n+1}p},y_{n+1})$.
\end{theorem}

\noindent\textbf{Proof.~}
The proof is an exercice, and the approach is similar to the proof
of Theorem 1.
$\Box$\medskip


\begin{theorem}
Let $(x_{11},\dots,x_{1p},y_1),\dots,(x_{n1},\dots,x_{np},y_n)$ 
be $n$ points in 
$\mathbb R^{p+1}$, with $n>p+1$. Let $L(k)$ and $O(k)$, for 
$k=1,\dots,n$, defined as in Section~\ref{sectionone}.
Then $L(k)+O(k)\le n-1$.
\end{theorem}

\noindent\textbf{Proof.~}
This is a consequence of the fact that for each of the $n$ subsets 
the scores are
shared among distinct points, and, by the above assumptions, 
a point cannot collect
a score for both $L$ and $O$, since, for $L$ it must have 
zero residual, and for $O$
a strictly positive absolute residual. And each point appears 
exactly $n-1$ times in the
$n$ subsets.
$\Box$\medskip


\section{The Algorithms}\label{algorithms}
In this section we propose two algorithms based on the previous
section. The aim of Algorithm 1 and Algorithm 2 is to detect 
leverage points and outliers respectively.


\bigskip
\noindent
\textbf{Algorithm 1.}
\textit{
\begin{enumerate}
\item Consider a data set $S$ of size $n$. Let $A$ and $B$ be
empty sets.
\item
Let $m$ be the size of $S$.
\item
Consider all $m$ subsets of $S$ of size $m-1$ and fit the
LAD-regression model for each subset.
\item
Compute $L(k)$ for each point of $S$ and select $k_1\in S$ which
maximizes $L(k)$.
\item
If $L(k_1)\ge\frac89(m-1)$ and $L(k_1)\ge\frac34(n-1)$ then move
$k_1$ into $B$ and move the eventual points of $A$ into $S$; 
otherwise move $k_1$ into $A$. 
\item
If the size of $S$ does not exceed $\frac9{10}n$ then the process
stops and $B$ represents the leverage points; otherwise go to step
2.
\end{enumerate}
}

\medskip

In this algorithm, the elements of $S$ are transferred in a set
$B$ of leverage points, or in a temporary set $A$ where points
that did not reached a sufficient score to be classified as
leverage points, are suitable to be reconsidered after that
another point has been detected. This trick avoids the {\it
masking effect}.

Discriminating values $\frac89(m-1)$,
$\frac34(n-1)$ and $\frac9{10}n$ for the score function $L$, 
have been empirically
determined, by testing on several data sets and several combinations
of values. They have, however, a natural interpretation. When
there is a unique leverage point, almost all $m$ regression models
detect it, so its $L$ score is near to the maximum. When there are
more leverage points, scores may be very different, and the {\it
masking effect} can produce relatively small scores. Finally, we
keep into account the size of $S$, to determine how many leverage
points a data set may have. The process stops when the size of set
$S$ does not exceed $\frac9{10}n$, so with our method a data set
cannot have more than $\frac1{10}n$ leverage points.


\bigskip
\noindent
\textbf{Algorithm 2.}
\textit{
\begin{enumerate}
\item
Consider a data set $S$ of size $n$. Let $C$ and $D$ be empty
sets. 
\item
Initialize the last maximum score (LMS) by $0$.
\item
Let $m$ be the size of $S$.
\item
Consider all $m$ subsets of $S$ of size $m-1$ and fit the
LAD-regression model for each subset.
\item
Compute $O(k)$ for each point of $S$ and select $k_1$ which
maximizes $O(k)$.
\item
If $O(k_1) = m-1$
\\
a) Then if $O(k_1)=LMS-1$ or $LMS=0$ then move $k_1$ into $D$,
put $LMS=O(k_1)$ and move the eventual points of $C$ into $S$;
otherwise the process stops and $D$ represents the outliers.
\\
b) Otherwise move $k_1$ into $C$.
\item
If the size of $S$ does not exceed $\frac45n$ then the process
stops and $D$ represents the outliers; otherwise go to step 3.
\end{enumerate}
}
\medskip

In this algorithm, the set $D$ contains the points classified as
outliers and the set $C$ contains the points that did
not reached the score to be classified as outliers, and that are suitable
to be reconsidered in further steps until the algorithm stops.

We can note the two proposed algorithms have a similar structure.
However, the main difference is a feature of Algorithm 2: the outliers
have decreasing scores $O_1$, $O_1-1$, $O_1-2$, and so on. The algorithm
stops when this sequence cannot be continued.



The process also stops when the size of set $S$ does not exceed
$\frac{4}{5}n$, so here a data set cannot have more than
$\frac{1}{5}n$ outliers.

Discriminating values for the score function $O$, have been also
empirically determined.

\section{Some Examples}\label{examples}

In this section we illustrate the proposed algorithms and compare
them with other two methods using several real and simulated data
sets.

One of the method is the P-R plot proposed by Hadi~\cite{hadi} to aid
in classifying observations as leverage points, outliers or
combinations of both. Some authors suggested that points with
$h_{ii}>\frac{2(p+1)}{n}$, where $h_{ii}$ is the {\it i}th
diagonal element of matrix $H$, 
$p$ is the number of predictors and
$n$ the number of observations, can be classified as leverage points
and the points with $\frac{r_i}{\hat{\sigma} \sqrt{1-h_{ii}}}>2$,
where $r_i$ is the residual of the {\it i}th observation, $h_{ii}$
is the {\it i}th diagonal element of matrix $H$ and ${\hat{\sigma}}$
is an estimator of standard deviation of the errors, can be classified
as outliers. In what follows the use of these 
suggested cut-off points to classify the
observations will be intended as \textit{classical methods}.

The first data set `Telephone' relate the number of international
telephone calls from Belgium (in tens of millions in minutes) to
the variable year for 24 years and can be found in ~\cite{rousseeuwleroy}. 
Cases 15-20 are unusually high and they are
outliers. The second one 'Hawkins' consists of 75 observations in
four dimensions, one response variable and three predictor
variables, and can be found in ~\cite{hawketal984}.
It has been constructed for the study of special pathological
phenomena in detection of outliers and leverage points and the
cases 1-10 are outliers and leverage points. The data set
'Scottish' describes how the record times (in seconds) of 35
Scottish Hill races is related to two predictor variables,
distance of race (in miles) and climb (in feet), and can be found
in ~\cite{hadi}. The data contain two clear outliers (observation 7
and 18). The last two data sets have been created by the authors. The
data set 'twovariables' consists of 56 observations on one
predictor variable and a response variable. The predictor was
created as uniform $(0,10)$ and the response variable  to be
consistent with the model $Y=X_1+4+\varepsilon$ with $\varepsilon\sim
N(0,1)$. Three observations (51--53) have been conceived as
leverage points and three others (54--56) as outliers. The other
data set 'threevariables' is the three variable equivalent to the
preceding one (two predictor variables).

Computation has been performed with a computer code in Splus, and the
results are summarized in Table 1.


\begin{table}[ht]
\begin{center}
{
\begin{tabular}{|l|l|c|c|}\hline\hline
Data & Method & Leverages & Outliers \\ \hline
Telephone &Classical Method&- & 20  \\
          &Hadi's Method & -&19, 20  \\
          & Our Results & -& 17-20  \\ \hline
Hawkins   & Classical Method  & 12-14 & 7, 11-14\\
          & Hadi's Method  &  14 & 7, 11-14\\
          & Our results & 3-6, 9, 10, 13 & 11-14 \\ \hline
Scottish  & Classical Method & 7, 11, 33, 35& 7, 18 \\
          & Hadi's Method  & 7, 11& 7, 18\\
          & Our Results & 11, 17, 35& 7, 18, 33\\ \hline
Twovariables & Classical Method  & 51-53 & 54-56\\
             & Hadi's Method & 51-53& 54-56 \\
             & Our results  & 52 & 54-56\\ \hline
Threevariables & Classical Method  & 18, 51-53& 54-56 \\
               & Hadi's Method  & 51-53 & 54-56\\
               & Our Results   & 51-53 & 9, 37, 54-56\\ \hline\hline
\end{tabular}
}
\mbox{  }\\\medskip
Table 1. \textsl{Detection of outliers and leverage points according to
different methods.}
\end{center}
\end{table}

\vspace{0.5 cm}
As we can see in Table 1, our proposed
method performed very well in detecting all outliers in the data
set 'Telephone'. The other methods failed to identify all of them
because the observation 19 and 20 mask all the others.

In the case of 'Hawkins' data, our proposed method as
well as the other two methods failed to identify the outliers. 
The outliers are all also swamped in the good cases 11--14. 
On the other hand,
our method detected almost all leverage points.

In the data set 'Scottish', Table 1 shows that all three methods
identified correctly the observation 7 and 18 as outliers. These
observations mask the observation 33 detected by our proposed
method. The observation 11 is suitably detected as leverage point
by all three methods, but there are others observations identified
as leverage points only in one or two methods.

For the simulated data sets, all methods performed very well in
detecting all outliers. However, our proposed method failed to
identify all leverage points.

\section{Conclusion}
The computation of the scores requires the determination 
of a certain number of LAD regression
models, and this is computationally longer than usual methods. 
However it is important to note that the 
principle is very simple, and takes into account natural requirements
for robustness in the detection of influential observations. 
Nowadays, the performances of 
common notebooks are largely
sufficient to perform in a few seconds the computations for 
the above examples, so the new tools are suitable for applications in
the stattistical methodology.

A code has been implemented in Splus and is available at the web site
of the second author
http://www.unine.ch/statistics/melfi/lad.html. A variety of data sets, 
including simulated  datasets used in Section~\ref{examples}, is also
available on the same web site.



\end{document}